\documentclass[10pt]{article}
\usepackage[utf8]{inputenc}
\usepackage[top=1in,bottom =1in, left=1in,right=1in]{geometry}

\usepackage{amsmath,amssymb,mathtools,bm,isomath,amsthm}
\usepackage{graphicx}
\usepackage{epstopdf}
\usepackage{color}
\usepackage[dvipsnames]{xcolor}
\usepackage{soul}
\usepackage{enumerate}
\usepackage{nicefrac}
\usepackage{enumitem}
\usepackage{mathrsfs}
\usepackage{scalerel,stackengine}
\usepackage{afterpage}
\usepackage{url}
\usepackage[affil-it]{authblk} 
\usepackage[toc,page]{appendix}
\usepackage{lineno}
\usepackage{physics}

\stackMath
\newcommand\reallywidehat[1]{%
	\savestack{\tmpbox}{\stretchto{%
			\scaleto{%
				\scalerel*[\widthof{\ensuremath{#1}}]{\kern-.6pt\bigwedge\kern-.6pt}%
				{\rule[-\textheight/2]{1ex}{\textheight}}
			}{\textheight}%
		}{0.5ex}}%
	\stackon[1pt]{#1}{\tmpbox}%
}

\newcommand{\Sc}{\mathcal{S}}
\newcommand{\Uc}{\mathcal{U}}
\newcommand{\bx}{{\boldsymbol x}}

\newcommand{\bb}{{\boldsymbol b}}
\newcommand{\bs}{{\boldsymbol s}}
\newcommand{\bz}{{\boldsymbol z}}

\begin{document}
	\setlength{\baselineskip}{15pt}
	
	\title{Operator approximation of the wave equation\\ based on deep learning of Green’s function}
	
	\author[1]{Ziad Aldirany}
	\author[2]{R\'egis Cottereau}
	\author[1]{Marc Laforest}
	\author[1]{Serge Prudhomme}
	\affil[1]{D\'epartement de math\'ematiques et de g\'enie industriel, Polytechnique Montr\'eal, Montr\'eal, Qu\'ebec, Canada}
	\affil[2]{Aix-Marseille Universit\'e, CNRS, Centrale Marseille, LMA UMR 7031, Marseille, France}
	
	\date{\today}
	
	\maketitle
	
	\begin{abstract}
	Deep operator networks (DeepONets) have demonstrated their capability of approximating nonlinear operators for initial- and boundary-value problems. One attractive feature of DeepONets is their versatility since they do not rely on prior knowledge  about the solution structure of a problem and can thus be directly applied to a large class of problems. However, convergence in identifying the parameters of the networks may sometimes be slow. In order to improve on DeepONets for approximating the wave equation, we introduce the Green operator networks (GreenONets), which use the representation of the exact solution to the homogeneous wave equation in term of the Green's function. The performance of GreenONets and DeepONets is compared on a series of numerical experiments for homogeneous and heterogeneous media in one and two dimensions.

	\end{abstract}
	
	
\section{Introduction}
\label{sec:introduction}

In the last few years, a large amount of work, see e.g.~\cite{jin2021nsfnets, bihlo2022physics, pettit2020physics, moseley2020solving}, has been devoted to using deep learning methods for the solution of PDE-based problems, such as in fluid dynamics, elasticity, meteorology, etc. These works have been motivated by the ability of deep neural networks to approximate large classes of functions in high dimension over complex domains~\cite{hornik1989multilayer, pinkus1999approximation}.
	
Prominent deep learning methods for solving partial differential equations rely on either learning the solution to the problem, as presented in~\cite{raissi2019physics, weinan2018deep, zang2020weak}, or learning the operators that describe the physical problem, as introduced in~\cite{lu2019deeponet, li2020fourier, li2020neural, wang2021learning}. In the first approach, the solution is approximated with a neural network by minimizing the residual of the PDE, as in so-called physics-informed neural networks (PINNs) introduced by~\cite{raissi2019physics}. The second approach learns the differential operator for a given family of parameters, e.g.\ deep operator networks (DeepONets) as presented in~\cite{lu2019deeponet}, thus allowing one to subsequently approximate the solution to a physical problem for a specific parameter in the vicinity of the trained parameters. The training of such a neural network can turn out to be quite expensive, but needs to be done only once. Computing the solution for a new parameter requires only one forward pass in the online phase, which is usually cost-effective. This makes the operator approximation method very attractive when the physical problem needs to be solved for a wide range of parameter values. For example, in seismology, uncertainties in the Earth's properties often require thousands of simulations to obtain those solutions that best fit the measured data.

In physics-based deep learning approaches, the training of the network does not require using any data as it is essentially based on the physics of the problem, in the sense that it approximates the solution to the partial differential equations along with the boundary and initial conditions. This is achieved by minimizing the residual associated with the partial differential equations. The most common approach is to consider a Least Squares approach, in which the loss function is defined in terms of the $L^2$-norm of the strong form of the residual~\cite{raissi2019physics, wang2021learning}. In case the solution of the problem lacks regularity, one can define the loss function by minimizing the energy potential in the case of symmetric problems, such as the Ritz formulation~\cite{weinan2018deep}, or by minimizing the norm of the residual functional, see e.g.~\cite{kharazmi2019variational, zang2020weak}. The choice of the loss function that one should consider for the solution of PDE problems with neural networks, is still viewed as an open question~\cite{wang20222}. However, evaluation of the residual has been made more amenable thanks to automatic differentiation~\cite{baydin2018automatic}. In this work, we will formulate the problem using the classical Least Squares method that involves the strong form of the residual due to its simplicity and convergence ability.

The main objective in this paper focuses on approximating the operator of the wave equation for a family of initial conditions. The architecture presented in DeepONets is general and can be applied to a large class of parametric PDEs. In order to improve on DeepONets, we propose here an approach based on the representation of the exact solution to the homogeneous wave equation on unbounded domains in terms of the Green's function (for details on Green's function, we refer the reader to~\cite{duffy2015green}). The method will be heretofore referred to as the Green operator networks (GreenONets). GreenONets provide an approximation of the operator of the wave equation in terms of the corresponding Green's function. Similar techniques based on the approximation of the Green's function have been recently considered for the solution of linear and non-linear operators, see e.g.~\cite{boulle2022data, gin2021deepgreen, li2020neural}. The architecture used in~\cite{boulle2022data} is similar to the one presented in our work. However, in their work, the authors approximate the operator for a family of forcing terms and the training is performed by minimizing the error between the approximated solution and the exact solution, which makes their approach data-based, while the training in our work is based on the physics of the problem. In~\cite{gin2021deepgreen}, a dual-autoencoder architecture is presented to approximate the operator for non-linear boundary value problems, by linearizing the problem and approximating the corresponding Green's function. The authors in~\cite{li2020neural} introduced the graph neural operator, which is inspired by the Green's function. However, the graph neural operator does not compute the Green's function but aims at learning a corresponding kernel function using an iterative architecture. We will illustrate on a series of numerical examples involving the homogeneous and heterogeneous wave equations that the approximation of the operator using GreenONets exhibits in general better results in terms of accuracy and convergence when compared to DeepONets.

The paper is organized as follows. We introduce in Section~\ref{sec:preliminaries} the model problem and some preliminary notations. In Section~\ref{sec:deeponets}, we briefly present DeepONets and its architecture. We then describe in Section~\ref{sec:greenonets} the main features of GreenONets. We compare the performance of GreenONets and DeepONets on several numerical results for the homogeneous and heterogeneous wave equation in Section~\ref{sec:numresults}. Finally, our main conclusions and potential extensions of the current research work are summarized in Section~\ref{sec:conclusion}.
    
\section{Preliminaries}
\label{sec:preliminaries}

We introduce here some preliminaries and notations in order to describe the notion of operator of the wave equation using neural networks. We first present the model problem and continue with a brief account of neural networks and the use of PINNs to solve initial boundary-value problems.
    
\subsection{Model problem}
\label{subsec:problem}
	
The linear wave equation describes small perturbations from the steady state of a system that locally behaves like an elastic body. The material system is entirely characterized in its domain $\Omega \subset \mathcal{R}^d$, $d=1,2$, or $3$, by the distribution of a bounded function $c(\bx)$ in $\Omega$. The perturbations are introduced either as initial displacements $u_0$, initial velocities $u_1$, or as perturbations entering the domain through its boundary $\partial \Omega$ at different times. For the sake of simplicity, we will assume that the boundary conditions are given by homogeneous Dirichlet boundary conditions on the displacement.
To be more specific, given a wave speed \(c(\bx)\), initial displacement \(u_0(\bx)\), initial velocity \(u_1(\bx)\), and a final time \(T>0\), the problem is to find the perturbation \(u(\bx,t)\), for all \(\bx\in\bar\Omega\) and \(t\in(0,T)\) such that
\begin{equation}
\partial_{tt} u(\bx,t) - c^2(\bx) \laplacian u(\bx,t) = 0, 
\quad \forall (\bx,t) \in \Omega \times (0,T),
\end{equation}
subjected to the initial and boundary conditions 
\begin{align}
u(\bx,0) = u_0(\bx), &\quad \forall \bx \in \Omega,\\
\partial_t u(\bx,0) = u_1(\bx), &\quad \forall \bx \in \Omega,\\
u(\bx,t) = 0, &\quad \forall (\bx,t) \in \partial \Omega \times (0,T).
\end{align}
Our goal is to obtain a neural network approximation of the (inverse) operator of the wave equation that would provide the solution \(u=u(\bx,t;u_0)\) for a family of initial conditions \(u_0\). In this case, the networks will be trained on a family on initial conditions generated by Gaussian random fields (GRF)~\cite{rasmussen2006gaussian}. The resulting neural network solution will allow one to compute an approximation to the wave equation for any initial condition. The hope is that this approximation should be accurate for initial conditions that are close to those used in the training phase. For the sake of simplicity, we shall focus mostly on the case where \(u_1 = 0\) but we will indicate how our approach can be extended to non-zero initial velocities $u_1$.

\subsection{Green's functions}

The Green's function of the operator defined previously is defined as the solution of the same problem with an impulse (localized in space) as initial condition. More precisely, 
\(g(\bx,t,\xi)\) is defined,
for all \(\bx\in\Omega\) and \(t\in(0,T)\), as the solution of
\begin{equation}
\partial_{tt} g(\bx,t,\xi) - c^2(\bx) \laplacian g(\bx,t,\xi) = 0, 
\quad \forall (\bx,t) \in \Omega \times (0,T),
\end{equation}
subjected to the initial and boundary conditions 
\begin{align}
\label{eq:IC0_green}
g(\bx,0,\xi) = \delta(\bx-\xi), &\quad \forall \bx \in \Omega,\\
\label{eq:IC1_green}
\partial_t g(\bx,0,\xi) = 0, &\quad \forall \bx \in \Omega,\\
g(\bx,t, \xi) = 0, &\quad \forall (\bx,t) \in \partial \Omega \times (0,T).
\end{align}

An interesting feature of the Green's function is that the solution $u(x,t)$ of the problem presented in the previous section can be obtained as a simple convolution with the initial condition $u_0(\bx)$:
\begin{equation}
\label{eqn:greensfct}
u(x,t) = \int_\Omega g(x,t,\xi) u_0(\xi) d\xi.
\end{equation}
Although we have considered here the case when $u_1(x)$ vanishes, an additional Green's function can be defined with the Dirac delta function on Eq.~\eqref{eq:IC1_green} rather than on Eq.~\eqref{eq:IC0_green}, so that the solution would become a sum of two convolutions by the principle of superposition.

A Green's function can be defined for any geometry of the domain and any distribution of properties $c(x)$. However, in the more simple case of an unbounded domain $\Omega=\mathbb{R}$ and homogeneous properties $c(x)=c$, the Green's function takes the simple form
\begin{equation}
g(x,t,\xi) = \dfrac{1}{2} \delta\big(ct - |x - \xi|\big).
\end{equation}
More examples of Green's function (for higher dimensions and bounded domains) can be found in~\cite{duffy2015green} or~\cite{Kausel2006green}.

\subsection{Neural networks}

Neural networks have been the subject of intensive research in the past decades~\cite{krizhevsky2017imagenet, goodfellow2020generative} and more recently have been used as a discretization approach for solving differential equations~\cite{sirignano2018dgm, raissi2019physics}.
By definition, a neural network maps an input into an output by a composition of linear and nonlinear functions, with adjustable weights and biases. The objective is usually to train the network by adjusting its weights and biases in order to minimize some measure of error between the output and the corresponding target values over a specific  training set. In this sense, the optimal neural network is very much like the least-squares fit of some fixed model to experimental data, but in contrast to a least-squares fit, the minimization problem might not always possess a unique solution. The resulting network can then be used as a predictive model that should hopefully provide accurate output when considering a wider set of input. There exist several neural network architectures, e.g. convolutional neural networks~\cite{krizhevsky2017imagenet}, feedforward neural networks~\cite{lecun2015deep}. We describe below the feedforward neural networks (FNN) that will be used later with DeepONets and GreenONets.
    
Let us consider a FNN with \(d\) hidden layers, each layer having a width \(N_i\), \(i=1,\ldots,d\), and let \(N_0\) denote the size of the input data and \(N_{d+1} \) the size of the output layer. Denoting the activation function by \(\sigma\), the neural network with input \( (\bx,t) \) and output \(u\) is defined as
\begin{equation}
\label{eqn:FNN}
\begin{aligned}
	&\text{Input layer: }&&  \bz_0 = (\bx,t), \\
	&\text{Hidden layers: }&&  \bz_{i} = \sigma ( W_i \bz_{i-1} + \bb_i), \quad i=1,\cdots,d, \\ 	
	&\text{Output layer: }&& u = W_{d+1} \bz_{d} + \bb_{d+1} ,
\end{aligned}
\end{equation}
where \(W_i\) is the weights matrix of size \(N_{i} \times N_{i-1}\) and \(\bb_i \) is the biases vector of size \( N_{i}\). For convenience, we will combine the weights and biases into the single parameter \( \theta \) of the neural network. In this work, we shall consider the \(\tanh\) activation function, but other activation functions could be used as well.

\subsection{Physics-informed neural networks and operator approximation}
     
We recall here the physics-informed neural network approach for solving partial differential equations, first introduced in~\cite{raissi2019physics}, as applied to the wave equation. We denote the residual associated with the partial differential equation of the wave equation as
\begin{equation}
\label{eqn:residual}
R\big(\bx,t,u\big) = \partial_{tt} u(\bx,t) - c^2(\bx) \laplacian u(\bx,t) , \quad \forall (\bx,t) \in \Omega\times(0,T),
\end{equation} 
introduce the residual associated with the Dirichlet boundary condition as:
\begin{equation}
\label{eqn:bc}
B\big(\bx,t,u\big) = u(\bx,t), \quad \forall (\bx,t)\in \partial \Omega\times(0,T), 
\end{equation} 
and the residuals associated with the initial conditions as:
\begin{align}
\label{eqn:ic}
&I_1\big(\bx,u\big) = u(\bx,t=0)-u_0(\bx), \quad \forall \bx\in\Omega, \\
\label{eqn:ic2}
&I_2\big(\bx,u\big) = \partial_t u(\bx,t=0), \quad \forall \bx\in\Omega. 
\end{align} 
We note that if \(u_1(\bx)\) is different from zero, then it should be subtracted from Eq.~\eqref{eqn:ic2}.

In PINNs, the goal is to obtain the solution \(u\) to the problem by approximating \(u\)  with a neural network \({u}_\theta(\bx)\). The training is usually performed by minimizing the following loss function:
\begin{equation}
\label{eq:losspinn}
L(\theta) 
= w_r \int_0^T\int_\Omega R(\bx,t,{u}_\theta)^2 dx dt
+ w_{bc} \int_0^T\int_{\partial \Omega} B(\bx,t,{u}_\theta)^2 dx dt+ w_{ic} \int_{\Omega} I_1(\bx,{u}_\theta)^2 + I_2(\bx,{u}_\theta)^2 dx,
\end{equation}
where \(w_r\), \(w_{bc}\), and \(w_{ic}\) are weighting coefficients.

The minimization is usually done using a gradient-based method, e.g.\ ADAM~\cite{DBLP:journals/corr/KingmaB14}, since the minimization problem is non-convex with respect to the trained parameters. Some of the main advantages of the PINNs is that it is a meshless method, and therefore we eliminate the process of mesh construction that can be very time consuming. Moreover, the implementation of the different types of boundary and initial conditions is similar in all cases and can be simply done by adding an extra weighted term in the loss function as presented in~\eqref{eq:losspinn}. 

\section{Deep operator networks}
\label{sec:deeponets}

If one seeks to calculate the solution using PINNs for several values of some parameters \(\bs\) of the model problem, one should have to recompute the approximate solution \(u_\theta\) for each instance of \(\bs\) by training the network from the beginning. This becomes quickly inefficient if the solution has to be evaluated for multiple values of the parameters, such as in a multi-query approach for uncertainty quantification or optimization. An alternative approach is to construct a surrogate model, in which the solution is searched as an operator acting on the parameters \(\bs\), i.e.\ \(u=Q(\bs)\). Similarly to PINNs, the operator \(Q\) can be approximated using deep learning by minimizing the loss function associated with the residuals of the partial differential equation and the boundary and initial conditions for a family of parameters. This approach becomes more attractive when the initial boundary value problem should be solved multiple times for different parameters, since it requires only one forward pass to compute the solution of a new parameter in the online phase.  

We briefly review the deep operator networks first introduced by~\cite{lu2019deeponet}. In this work, we will be learning the operator from the partial differential equation and the initial and boundary conditions. In other words, the physics-informed DeepONets, described in~\cite{wang2021learning}, will be presented.

\begin{figure} [tp!]
\centering
\includegraphics[width=0.8\linewidth]{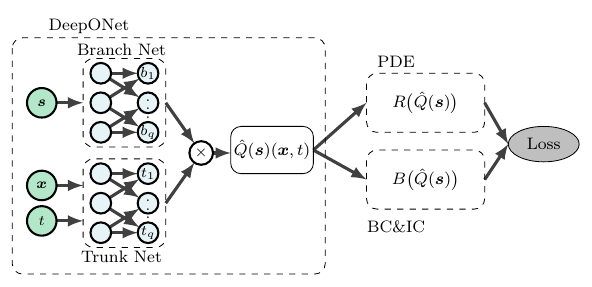}
\caption{Illustration of the architecture of the unstacked DeepONets. It consists of two networks, the branch net and the trunk net. The branch net takes as input the input vector while the trunk net takes the input coordinates. Their outputs are merged with dot product to give the approximated operator. The network is trained to minimize the loss function that consists of the residuals of the partial differential equation and the initial and boundary conditions.}
\label{fig:deepOnet}
\end{figure}

Given Banach spaces \(\Uc\) and \(\Sc\), we want to learn the operator \(Q : \Sc \to \Uc \) such that for any input parameter \( s \in \Sc\) (which in our case represents an initial condition), \(Q(s) \equiv u \in \Uc\) is the solution of problem~\eqref{eqn:residual} with the boundary conditions~\eqref{eqn:bc}. In order to approximate the operator \(Q\), we present the unstacked DeepONets architecture originally introduced in~\cite{lu2019deeponet} and schematically shown in Figure~\ref{fig:deepOnet}. We start by defining the input vector \(\bs\) of initial conditions \( [u_0(\bx_i)]_{i=1,\ldots,m}\) evaluated at a collection of $m$ points \({\{\bx_i\}_{i=1}^m}\), known as sensors. Then, as illustrated in Figure~\ref{fig:deepOnet}, the operator is approximated as
\begin{equation}
\hat{Q}(\bs)(\bx, t) = \sum_{k=1}^q b_k\big(\bs\big)t_k(\bx, t),
\end{equation}
where \(\{b_k\}_{k=1}^q\) is the output of the branch network that takes \(\bs\) as an input and \(\{t_k\}_{k=1}^q\) is the output of the trunk network that takes \((\bx,t)\) as an input. The value of $q$ will be chosen in the numerical experiments as the width of the layers in the neural network. We will consider a simple FNN for both the branch and trunk networks. We note that a convolutional neural network~\cite{krizhevsky2017imagenet} could also be used for the branch network when working with uniformly distributed sensors. 

We consider here the physics-informed DeepONets, where the network is trained by penalizing the residuals associated with the governing partial differential equation and with the initial and boundary conditions for a family of $N$ input functions \(\{\bs^{(i)}\}_{i=1}^N\). The loss function in this case reads:
\begin{multline}
\label{eq:lossDeepONet}
L(\theta) 
= \dfrac{1}{N} \sum_{i=1}^N \bigg[ w_r \int_0^T\int_\Omega R\big(\bx,t,\hat{Q}(\bs^{(i)})(\bx, t)\big)^2 dx dt 
+ w_{bc} \int_0^T\int_{\partial \Omega} B\big(\bx,t,\hat{Q}(\bs^{(i)})(\bx, t)\big)^2 dx dt \\
+ w_{ic} \int_{\Omega} I_1\big(\bx,\hat{Q}(\bs^{(i)})(\bx, t)\big)^2 + I_2\big(\bx,\hat{Q}(\bs^{(i)})(\bx, t)\big)^2 dx \bigg].
\end{multline}
	
\section{Green operator networks}
\label{sec:greenonets}

The architecture presented in the DeepONets is a general architecture that works can accommodate different problems with different input parameters. Our objective here is to develop an approach that improves upon the efficiency of the DeepONets for certain types of parameters. We will focus on the solution of the wave equation for homogeneous and heterogeneous materials, as presented in Section~\ref{subsec:problem}. We thus propose the Green operator networks (GreenONets), that approximate the Green's function of the operator, to solve the aforementioned problem. 
	
\begin{figure}[t!]
\centering
\includegraphics[width=0.9\linewidth]{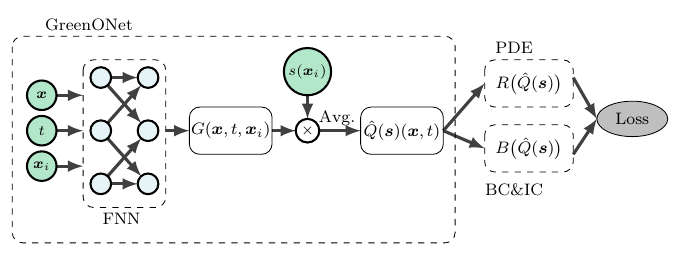}
\caption{Illustration of the architecture of the GreenONets. A FNN takes as input the coordinates and sensor points and outputs an approximated Green's function \(G\) with respect to each sensor point. Then, the operator is computed by averaging the product of the Green's function and the input function over the sensor points. The network is then trained to minimize the loss function that consists of the residuals associated with the partial differential equation and the initial and boundary conditions.}
\label{fig:greenOnet}
\end{figure}

In this work, we are interested in learning the operator of the wave equation for different initial conditions, i.e.\ the input function is defined as \(s = u_0\). Therefore, instead of using a general architecture as the DeepONets, we introduce the Green operator networks, shown in Figure~\ref{fig:greenOnet}, as a discrete approximation of the integral in Eq.~\eqref{eqn:greensfct}. The GreenONet is defined as
\begin{equation}
\label{eqn:gon}
\hat{Q}(\bs)(\bx, t) = \dfrac{1}{m}\sum_{i=1}^m G (\bx, t,\bx_i) u_0(\bx_i) ,
\end{equation}
where \(G\) is a simple feedforward neural network. We note that the formulation of GreenONets depends explicitly on the sensor points \({\{\bx_i\}_{i=1}^m}\). Similar to the physics-informed DeepONets, the GreenONets are trained by minimizing the loss function~\eqref{eq:lossDeepONet}. We notice that in our approximation, the solution is zero if the initial condition is zero. However, if the Dirichlet boundary condition or the initial speed are different from zero, we should add a term in~\eqref{eqn:gon}, that should be approximated by a new network, to compensate for these conditions.

Although the exact solution in~\eqref{eqn:greensfct} is only presented for a homogeneous material and an unbounded domain, the following numerical results show that GreenONets yield better results when compared to DeepONets for bounded domains with homogeneous or heterogeneous properties.
	
\section{Numerical results}
\label{sec:numresults}

In this section, we approximate the operator of the wave equation for homogeneous and heterogeneous materials in the case of a family of initial conditions, in order to show the effectiveness of GreenONets when compared to DeepONets.

We consider \(\Omega = (-1,1)^d\) and we set \(u_1(x) = 0\) in all cases. In order to define the training set, on which we want to minimize our residuals, we start by defining a family of $N$ input functions \(\{\bs^{(i)}\}_{i=1}^N\). For each \(\bs^{(i)}\), we randomly define \( \big\{\big(\bx_{ic,j}^{(i)},0\big) \big\}_{j=1}^{P_{ic}}\) on which we will penalize the initial conditions, \( \big\{ \big( \bx_{bc,j}^{(i)}, t_{bc,j}^{(i)}\big) \big\}_{j=1}^{P_{bc}}\) on the boundary of the domain at different times, for which we will penalize the boundary conditions, and \( \big\{\big(\bx_{r,j}^{(i)} ,t_{r,j}^{(i)}\big) \big\}_{j=1}^{P_r}\) on \((-1,1)^d\times (0,T)\), for which we will penalize the bulk residual. Using the sampling points to estimate numerically the integrals, the loss function~\eqref{eq:lossDeepONet} is approximated as 
\begin{equation}
\label{eqn:loss}
\mathcal L(\theta) = w_r \mathcal L_r(\theta)+w_{bc} \mathcal L_{bc}(\theta)+w_{ic} \mathcal L_{ic}(\theta),
\end{equation}
where 
\[
\begin{aligned}
&\mathcal L_r(\theta) = \dfrac{1}{NP_r} \sum_{i=1}^N \sum_{j=1}^{P_r}  \bigg|\partial_{tt} \hat{Q}\big(\bs^{(i)}\big) \big( \bx_{r,j}^{(i)}, t_{r,j}^{(i)}\big) - c\big(\bx_{r,j}^{(i)}\big)^2 \laplacian\hat{Q}\big(\bs^{(i)}\big)\big(\bx_{r,j}^{(i)}, t_{r,j}^{(i)}\big)\bigg)\bigg|^2, 
\\
&\mathcal L_{bc}(\theta) = \dfrac{1}{NP_{bc}} \sum_{i=1}^N \sum_{j=1}^{P_{bc}}  \bigg| \hat{Q}\big(\bs^{(i)}\big)\big(\bx_{bc,j}^{(i)}, t_{bc,j}^{(i)} \big) \bigg |^2, 
\\	
&\mathcal L_{ic}(\theta) = \dfrac{1}{NP_{ic}} \sum_{i=1}^N \sum_{j=1}^{P_{ic}}  \bigg| \hat{Q}\big(\bs^{(i)}\big)\big(\bx_{ic,j}^{(i)} , 0\big)  - u_0(\bx_{ic,j}^{(i)})\bigg |^2 + \bigg| \partial_t \hat{Q}\big(\bs^{(i)}\big)\big(\bx_{ic,j}^{(i)} ,0\big)\bigg |^2.
\end{aligned}
\]
Again, \(w_r\), \(w_{ic}\), and \(w_{bc}\) are the weighting coefficients. The initial conditions \(\bs^{(i)}\) are randomly sampled from a Gaussian random field (GRF), as presented by~\cite{lu2019deeponet}, with a defined length scale \(l\). In the following experiments, the FNNs in the DeepONets and GreenONets are defined with \(d = 6\) hidden layers and \(N_i = 50\) for all hidden layers. The loss function is minimized using the ADAM optimizer~\cite{DBLP:journals/corr/KingmaB14} with the default hyper-parameters, while considering different learning rates for each experiment.

\subsection{One-dimensional problems}
	
\begin{figure} [ht!]
\centering
\includegraphics[width=0.38\linewidth]{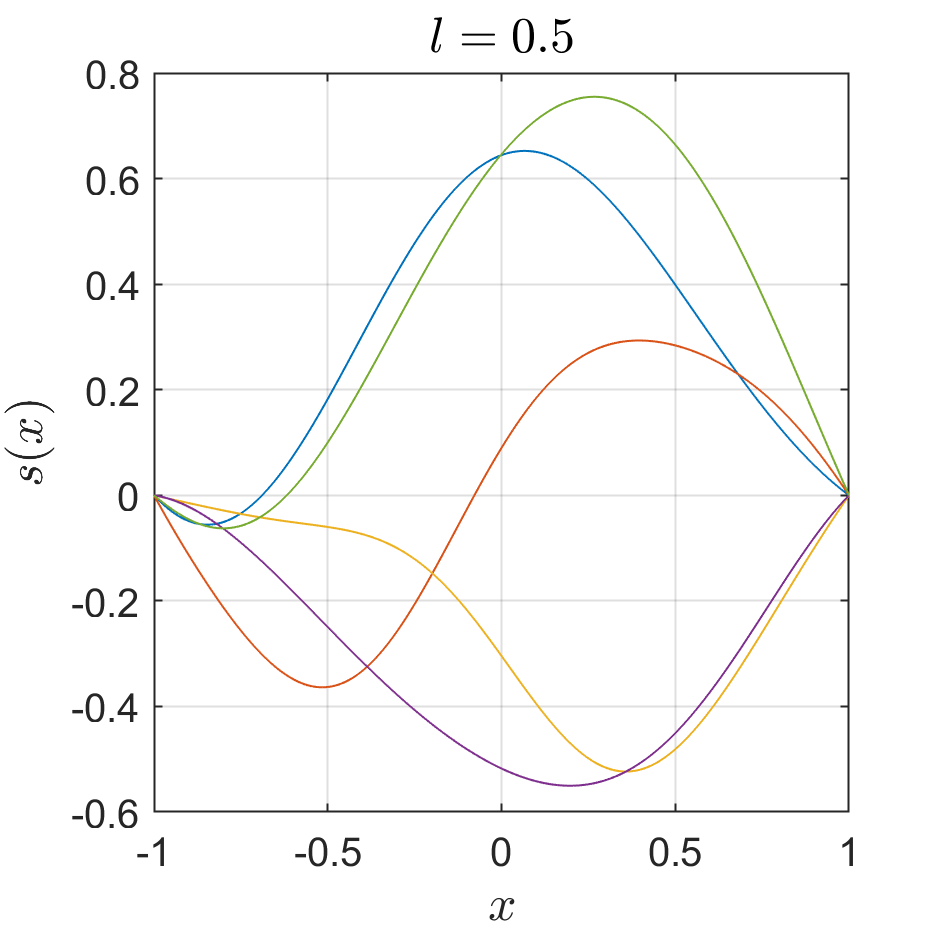}
\includegraphics[width=0.38\linewidth]{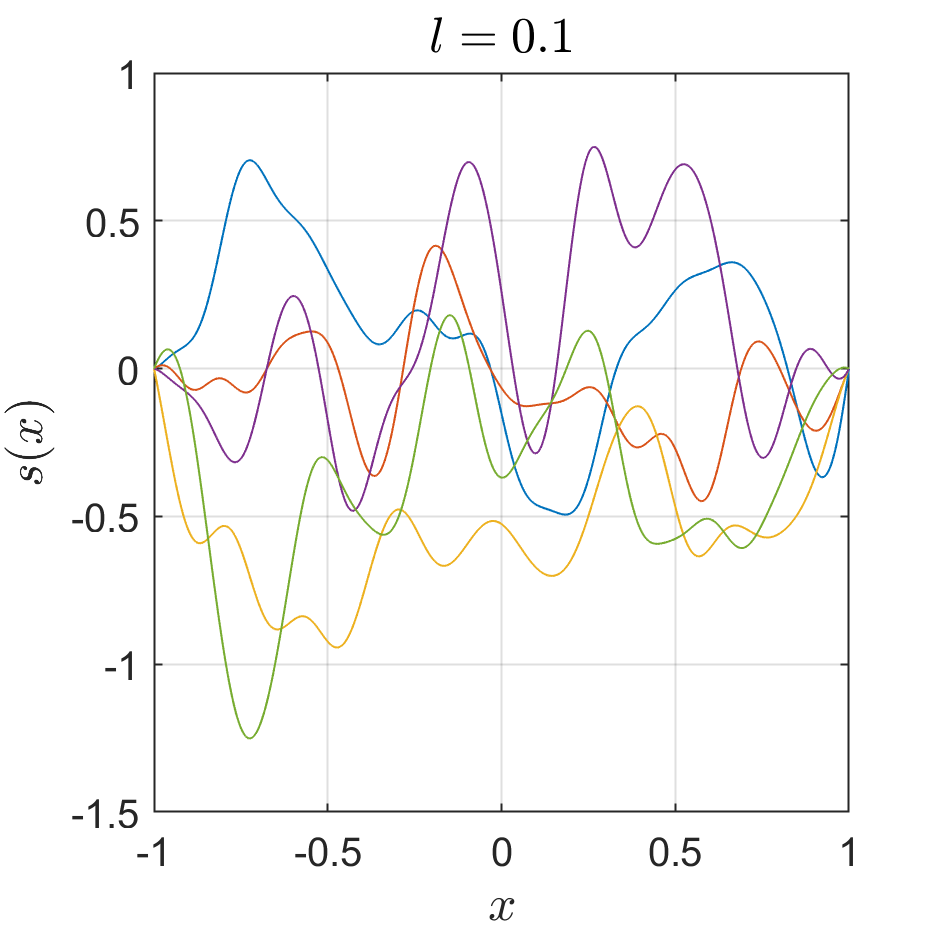}
\caption{Examples of the input functions \(s^{(i)}\) for a length scale \(l = 0.5\) (left) and \(l = 0.1\) (right).}
\label{fig:GRF}
\end{figure}

We start by comparing the GreenONets with the DeepONets for the one-dimensional case, i.e.\ \(d = 1\). The input functions \(s\) are defined by a GRF and then modified to verify the homogeneous Dirichlet boundary conditions by subtracting the proper linear function. Figure~\ref{fig:GRF} shows examples of the modified GRF for length scales $l=0.5$ and $l=0.1$. The numerical comparison is performed for homogeneous and heterogeneous material properties with different length scales in the GRF. 

\subsubsection{Homogeneous case with a length scale of 0.5}
\label{subsec:exp1}
	
We first approximate the operator for a homogeneous material using DeepONets and GreenONets. The sensor points are chosen uniformly with \(m=21\) while the input parameters \(\bs\) are generated using a GRF with length scale \(l=0.5\). We take the initial learning rate as \(10^{-3}\) and let it decrease with a rate of 0.9995 at each epoch. In this example, we choose \(N = 1000\) and \(P_r=P_{bc}=P_{ic}=10\). The weights in~\eqref{eqn:loss} are set to \(w_r = 0.1\) and \(w_{ic} = w_{bc} = 10\). The training is done for 5000 epochs with 16 mini-batches.
	
\begin{figure}[ht!]
\centering
\includegraphics[width=0.4\linewidth]{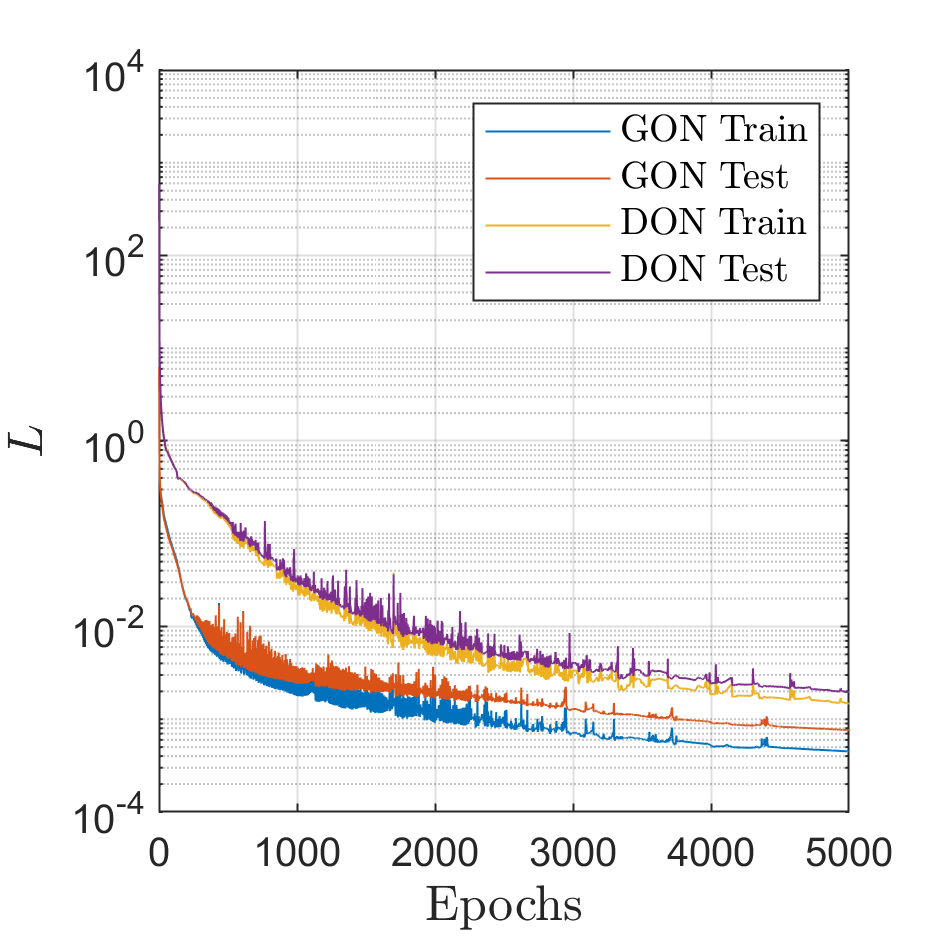}
\caption{The evolution of the loss function on the training and testing sets during the training with GreenONets and DeepONets for the example of Section~\ref{subsec:exp1}.}
\label{fig:loss1}
\end{figure}
	
\begin{figure} [ht!]
\centering
\includegraphics[width=0.32\linewidth]{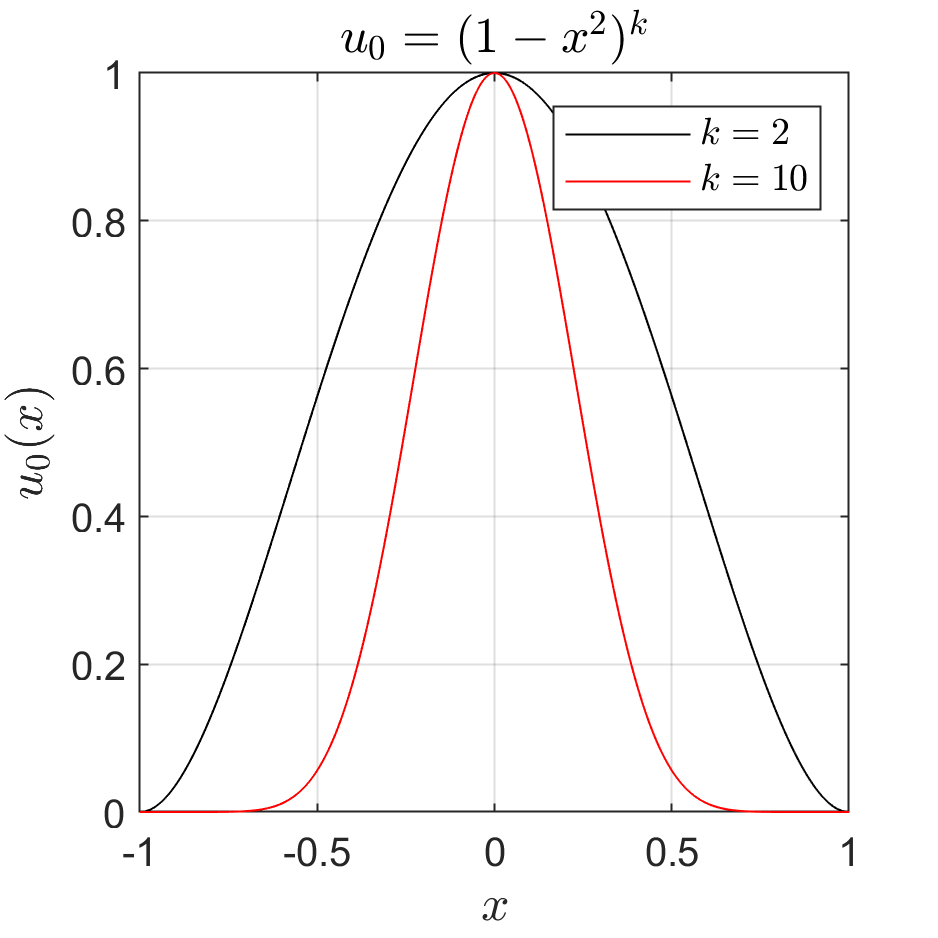}
\includegraphics[width=0.32\linewidth]{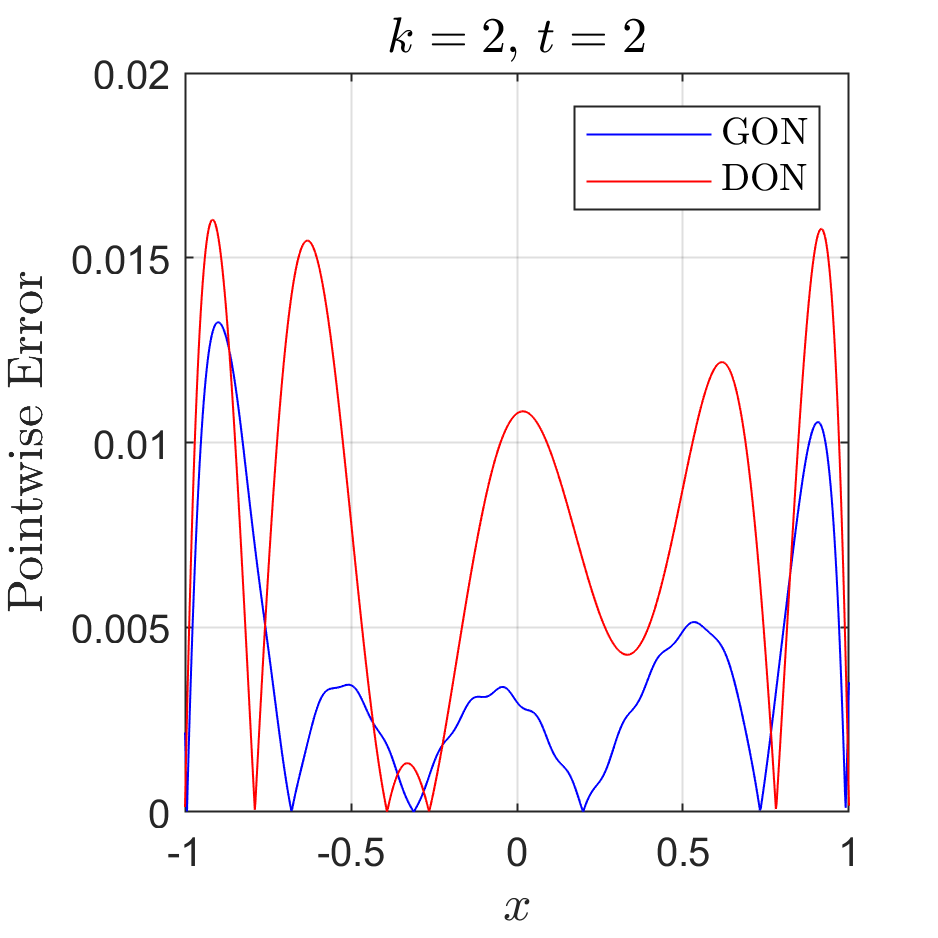}
\includegraphics[width=0.32\linewidth]{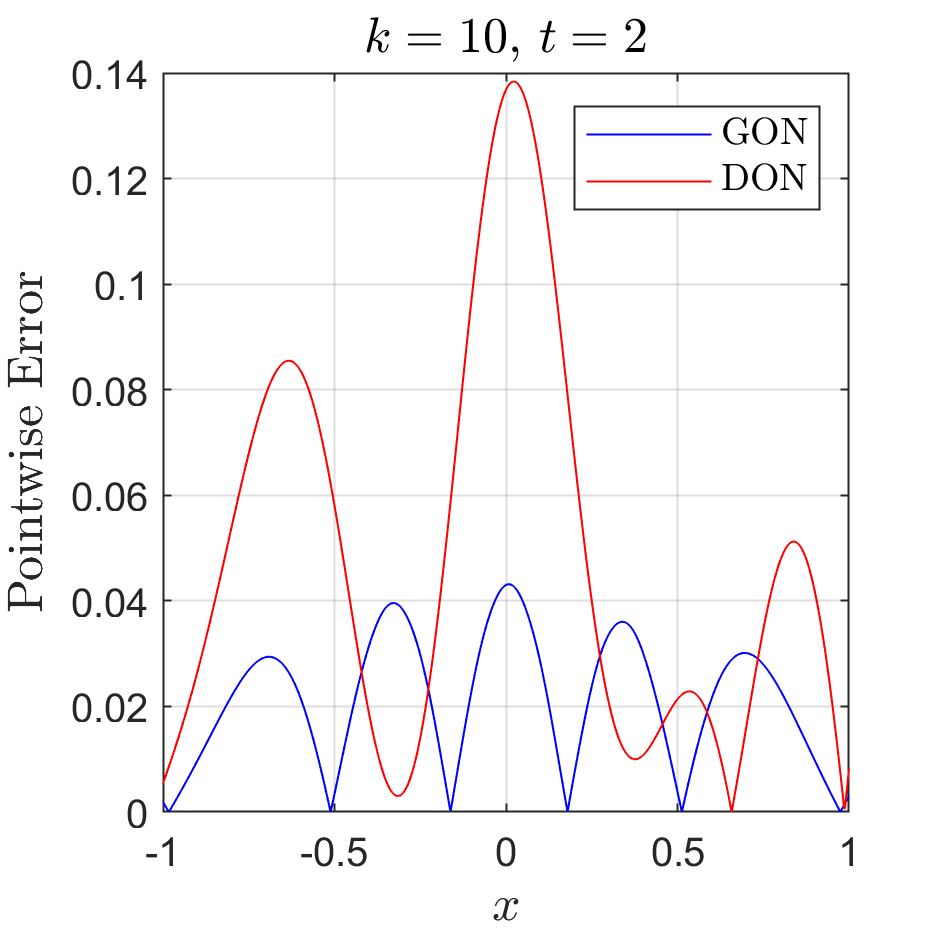}
\caption{Example of Section~\ref{subsec:exp1}: (left) Initial conditions with which we test the networks. (middle) Pointwise error at \(t=2\) for \(k=2\) using GreenONets and DeepONets. (right) Pointwise error at \(t=2\) for \(k=10\) using GreenONets and DeepONets.}
\label{fig:func_err1}
\end{figure}

Figure~\ref{fig:loss1} compares the evolution of the loss function on the training and testing sets for the GreenONets and DeepONets. We observe that the loss functions decrease faster with GreeONets and after 5000 epochs we have smaller losses when compared to the DeepONets loss functions. In order to verify our operators, we compute the solution at \(t=2\) for the initial conditions \(u_0(x) = (1-x^2)^k\), with \(k=2\) and \(k=10\), as shown in Figure~\ref{fig:func_err1} (left). We observe in Figure~\ref{fig:func_err1} (middle), that the pointwise error at \(t=2\) for \(k=2\) is slightly larger when using DeepONets. However, as shown in Figure~\ref{fig:func_err1} (right), for \(k=10\) the DeepONets solution exhibits a maximum pointwise error of 0.14 while the maximum pointwise error for the GreenONets solution is 0.04. Therefore, the GreenONets solutions seem to generalize better at higher frequencies.  
	
\subsubsection{Homogeneous case with a length scale of 0.1}
\label{subsec:exp2}
	
Here, we solve the same problem as in the previous section but the input parameters \(\bs\) are defined using a GRF with length scale \(l=0.1\). In other words, we now compare the two methods for higher frequency solutions. The sensor points are defined uniformly with \(m=60\). We initialize the learning rate to \(5\times10^{-4}\) and let it decrease with a rate of 0.999 at each epoch. In this example, we take \(N = 3000\), \(P_{r}=30\), and \(P_{bc}=P_{ic}=3\). The weights associated with each component of the loss function~\eqref{eqn:loss} are \(w_r = 0.2\) and \(w_{ic} = w_{bc} = 100\). We train both networks for 2000 epochs with 128 mini-batches.
	
\begin{figure} [ht!]
\centering
\includegraphics[width=0.32\linewidth]{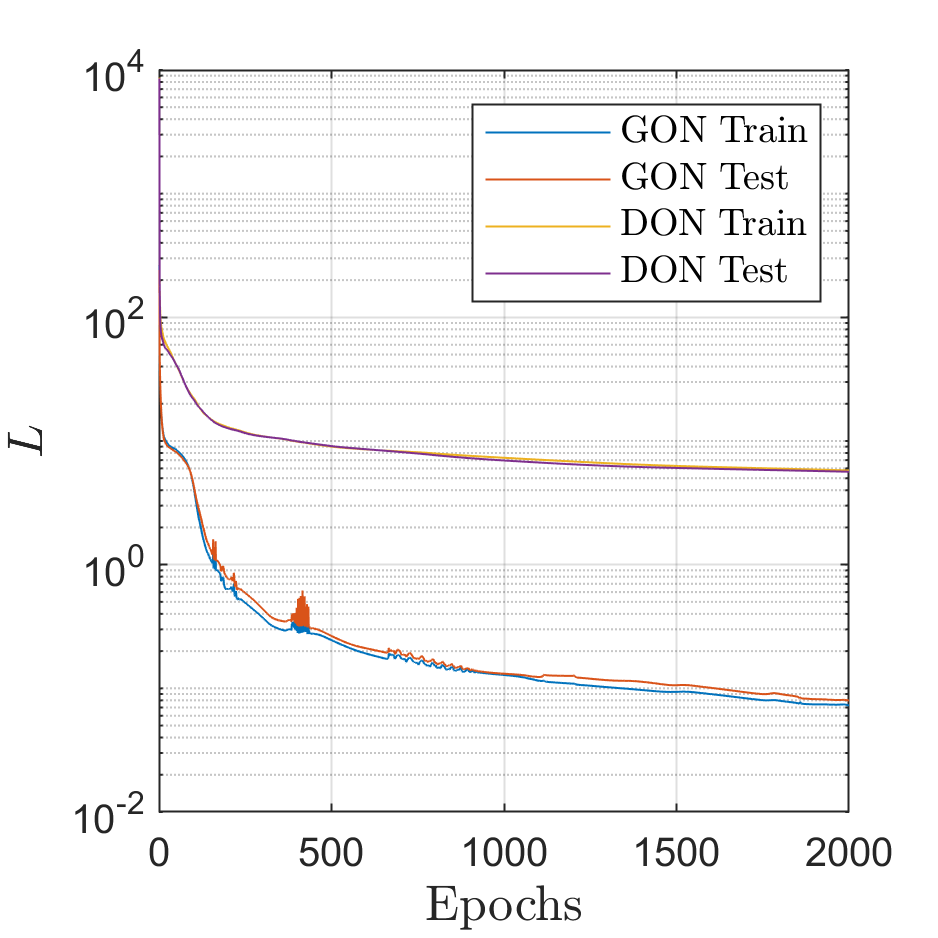}
\includegraphics[width=0.32\linewidth]{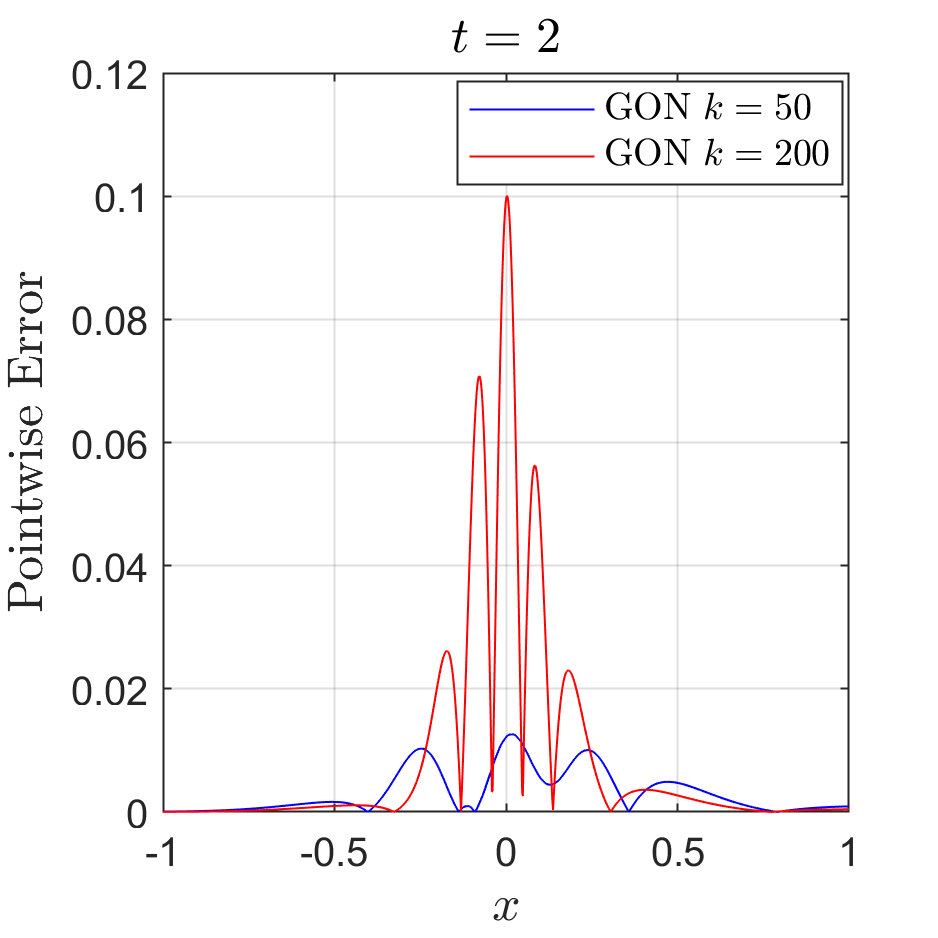}
\includegraphics[width=0.32\linewidth]{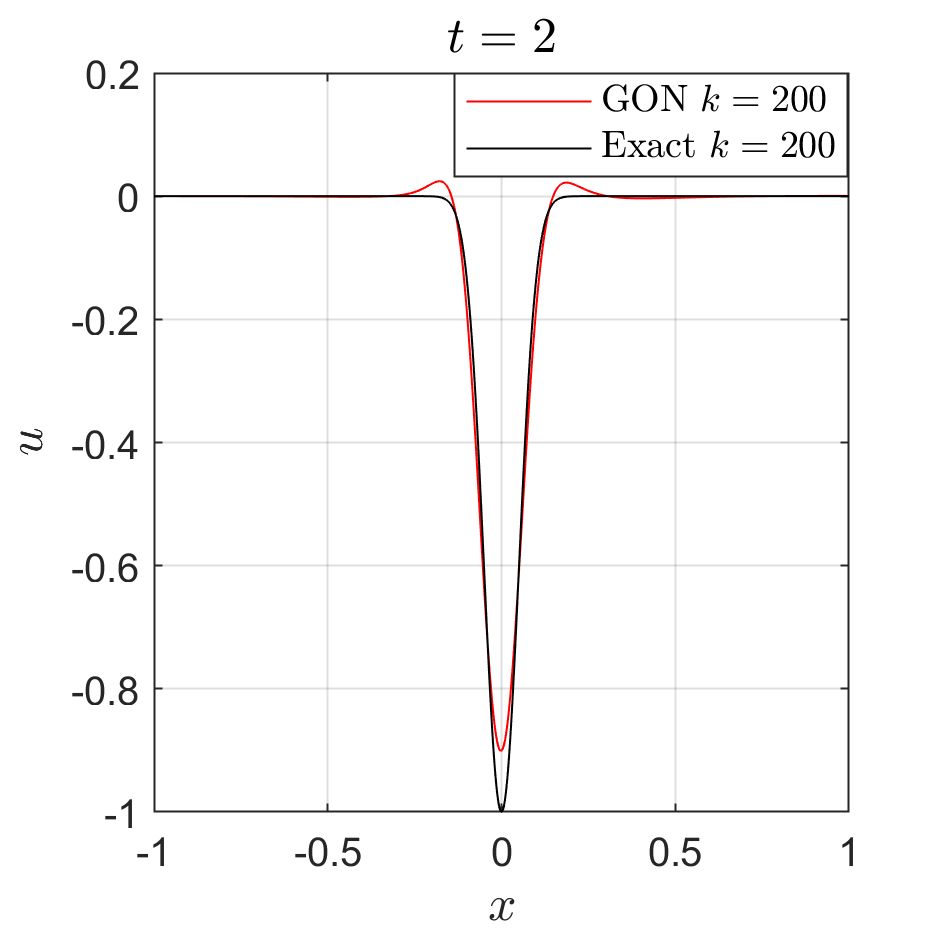}
\caption{Example of Section~\ref{subsec:exp2}: (left) Evolution of the loss function on the training and testing sets with GreenONets and DeepONets. (middle) Pointwise error at \(t=2\) for \(k=50\) and \(k=200\) using GreenONets. (right) Comparison of the solution obtained by GreenONets at \(t=2\) for \(k=200\) with the exact solution.}
\label{fig:loss2}
\end{figure}
	
We observe in Figure~\ref{fig:loss2} (left) that, using the same hyper-parameters, the loss functions for the GreenONets attain \(8\times10^{-2}\) in 2000 epochs while those for the DeepONets plateau earlier. In Figure~\ref{fig:loss2} (middle), we show the pointwise errors of the GreenONets solutions for \(u_0 = (1-x^2)^k\), with \(k=50\) and \(k=200\). The maximum pointwise error is around 0.01 for \(k=50\) and around 0.1 for \(k=200\). The pointwise error of the DeepONets solutions is not available since the loss functions did not converge. To better characterize the error for \(k=200\), we compare the solution using GreenONets to the exact solution at \(t=2\) in Figure~\ref{fig:loss2} (right). We remark that the large errors are close to the propagating wave and did not spread in the rest of the solution. 
	
\begin{figure} [ht!]
\centering
\includegraphics[width=0.4\linewidth]{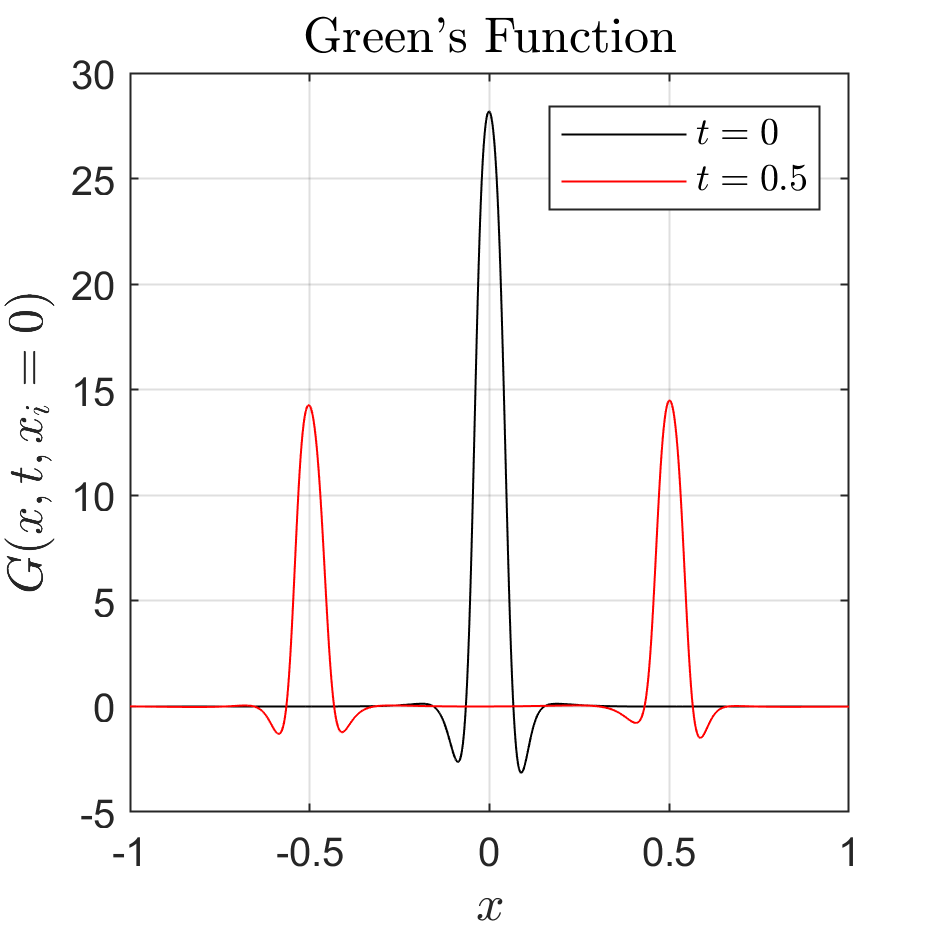}
\caption{Example of Section~\ref{subsec:exp2}: The approximated Green's function \(G(x,t,x_i) \), for \(x_i = 0\) at \(t = 0\) and \(t = 0.5\).}
\label{fig:Green2}
\end{figure}
	
We show in Figure~\ref{fig:Green2} the approximated Green's function \(G(x,t,x_i)\), as computed by Equation~\eqref{eqn:gon}, for \(x_i = 0\) at \(t = 0\) and \(t = 0.5\). As expected, we observe that \(G\) approximates a Dirac delta function around \(x = 0\) at \(t = 0\), which splits into two functions at \(t = 0.5\) with half the amplitude of the original one.
	
\subsubsection{Heterogeneous case with a length scale of 0.3}
\label{subsec:exp3}
	
In this section, we use the GreenONets to approximate the operator of the heterogeneous wave equation. We consider the wave speed \(c(x)^2 = 1 + \mathcal{H}(x-0.5)\), where \(\mathcal{H}\) is the Heaviside function. We consider a family of initial conditions with a length scale \(l=0.3\) and use uniformly distributed sensor points with \(m = 30\). The learning rate is \(10^{-3}\) and decreases with a rate of 0.9995 per epoch.  We set \(N = 2000\), \(P_r=15\), and \(P_{bc}=P_{ic}=3\). The weights of the loss functions are \(w_r = 1\) and \(w_{ic} = w_{bc} = 100\). We divide our training set into 32 mini-batches and train the networks for 2500 epochs.
	
\begin{figure} [ht!]
\centering
\includegraphics[width=0.3\linewidth]{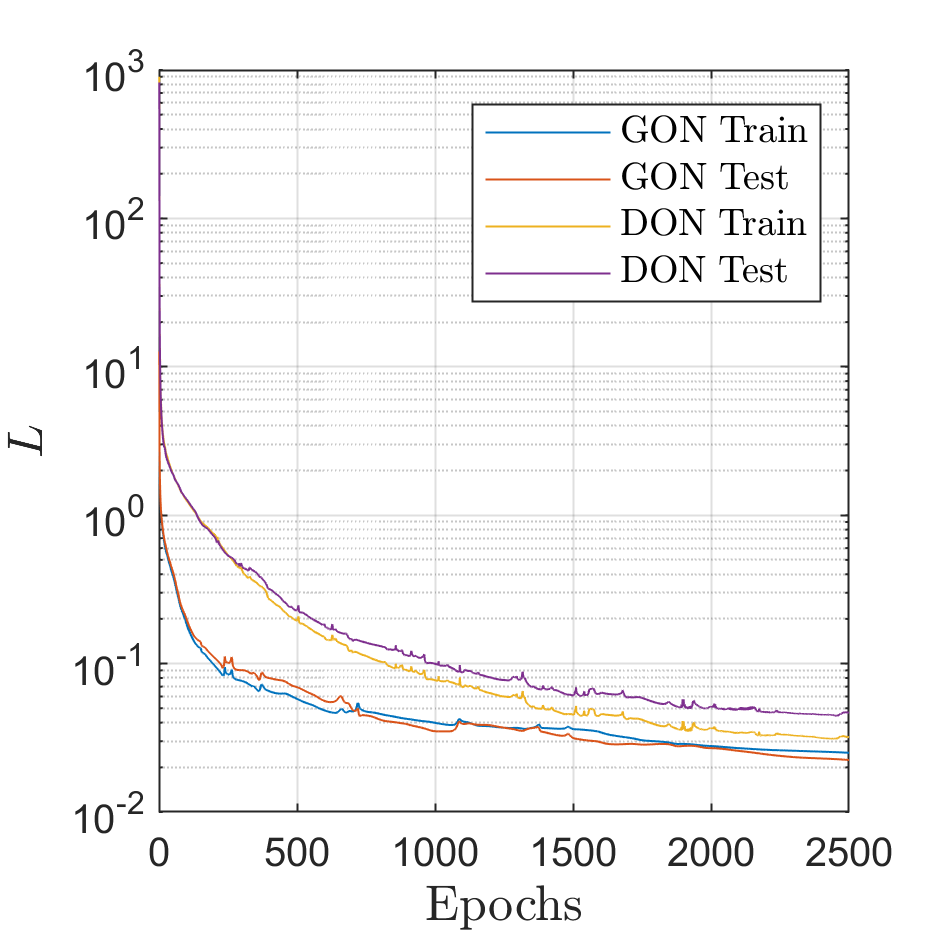}
\includegraphics[width=0.3\linewidth]{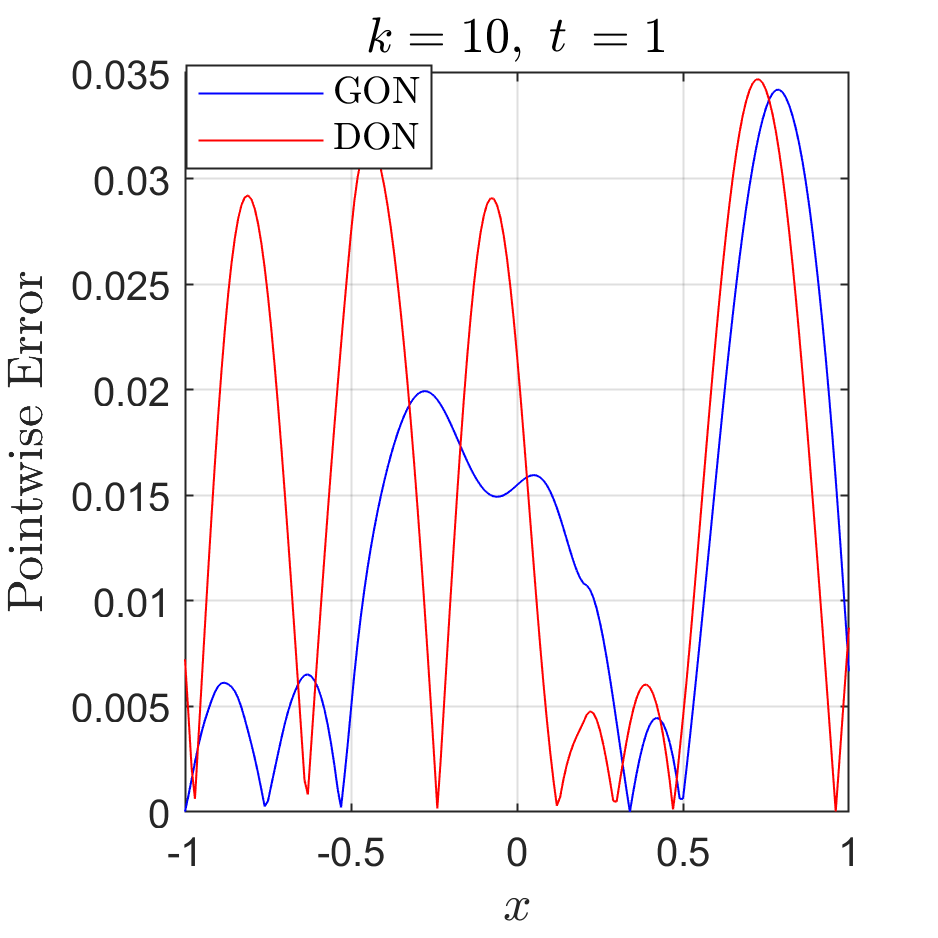}
\includegraphics[width=0.3\linewidth]{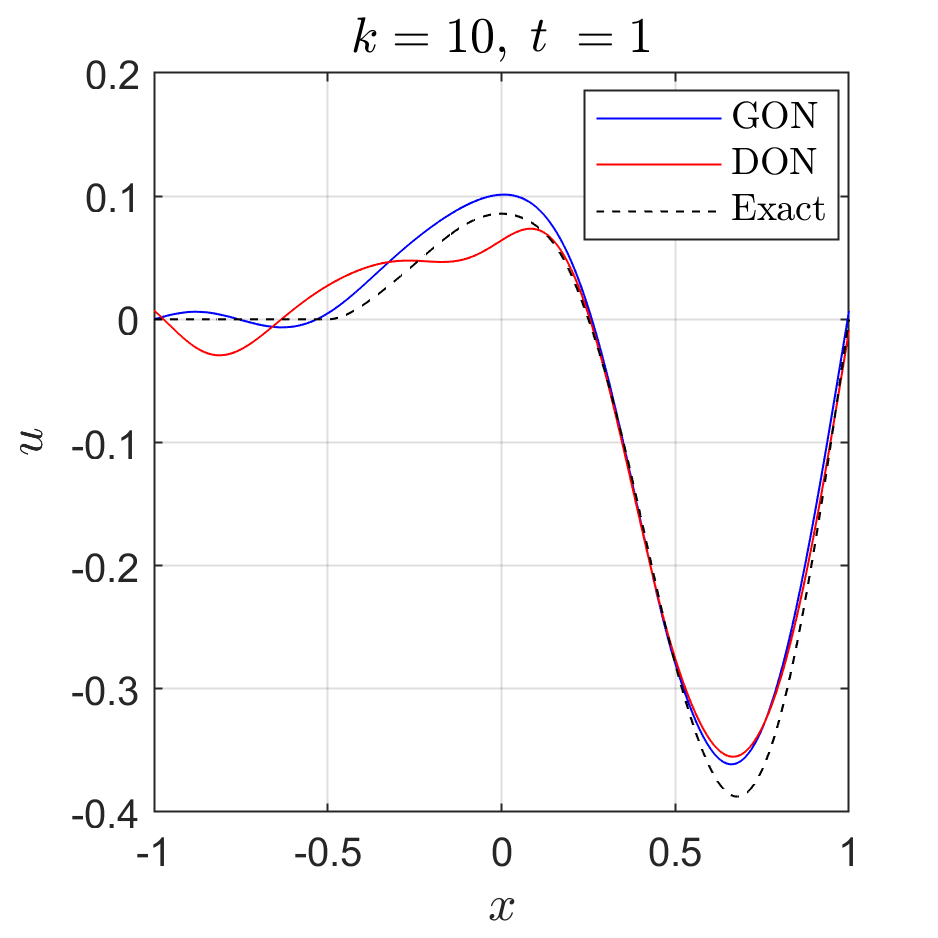}
\caption{Example of Section~\ref{subsec:exp3}: (left) Evolution of the loss function on the training and testing sets with GreenONets and DeepONets. (middle) Pointwise error at \(t=1\) for \(k=10\), using GreenONets and DeepONets. (right) Comparison of the solutions obtained by GreenONets and DeepONets at \(t=1\) for \(k=10\) with an overkill solution using the spectral element method.}
\label{fig:loss3}
\end{figure}
	
As shown in Figure~\ref{fig:loss3} (left), the loss functions on the training set and testing set decrease faster using GreenONets than with DeepONets. Figure~\ref{fig:loss3} (middle) shows the pointwise error in the solution at \(t=1\) using GreenONets and DeepONets with an initial condition \(u_0=(1-x^2)^{10}\). The maximum pointwise errors are similar for both networks with a value close to 0.035. We also plot the solutions at \(t=1\) along with the exact solution -- actually, an overkill solution using the spectral element method, see e.g.~\cite{Matlabcode,seriani1994,aldirany2022}, and references therein -- in Figure~\ref{fig:loss3} (right). We remark that the errors in the solution of the GreenONets remain localized around the main pulses and are proportional to the wave amplitude. However, in the case of DeepONets, the pointwise errors have the tendency to spread over the whole domain. Therefore, one could conclude that GreenONets tend to provide better approximations of the propagating waves without introducing large errors away from the main pulses. 
	
\subsection{A two-dimensional example}
\label{subsec:exp2d}
	
In this section, we present some numerical results obtained with DeepOnets and GreenOnets for the two-dimensional homogeneous wave equation with \(c(x,y) =1\), where \((x,y)\) denote the spatial coordinates. The input functions \(s\) are defined as
\[
s(x,y) = (1-x^2)(1-y^2)h(x,y),
\]
where \(h\) is randomly sampled from a zero-mean Gaussian random field with a length scale \(l = 1\). We choose \(m=49\) sensor points uniformly distributed in \(\Omega=[-1,1]\times[-1,1]\). We consider an initial learning rate \(10^{-3}\) that decreases with a rate 0.9995 at each epoch. The training set is defined with \(N = 50000\), and  \(P_{r} = P_{bc}=P_{ic} = 1\). The weights in Equation~\eqref{eqn:loss} are set to \(w_r = 1\) and \(w_{ic} = w_{bc} = 100\). The networks are trained for 400 epochs with 100 mini-batches.
	
\begin{figure}[ht]
\centering
\includegraphics[width=0.4\linewidth]{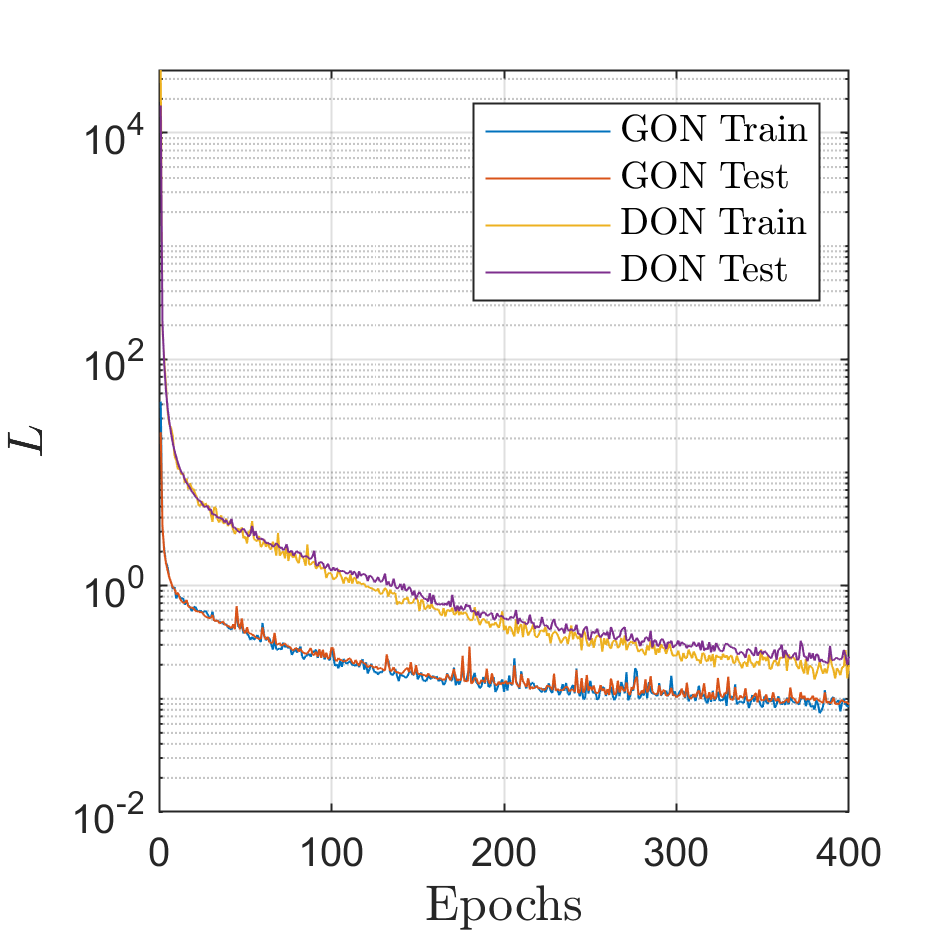}
\caption{Example of Section~\ref{subsec:exp2d}: Evolution of the loss function on the training and testing sets.}
\label{fig:loss2d}
\end{figure}

\begin{figure}[th]
\centering
\includegraphics[width=0.45\linewidth]{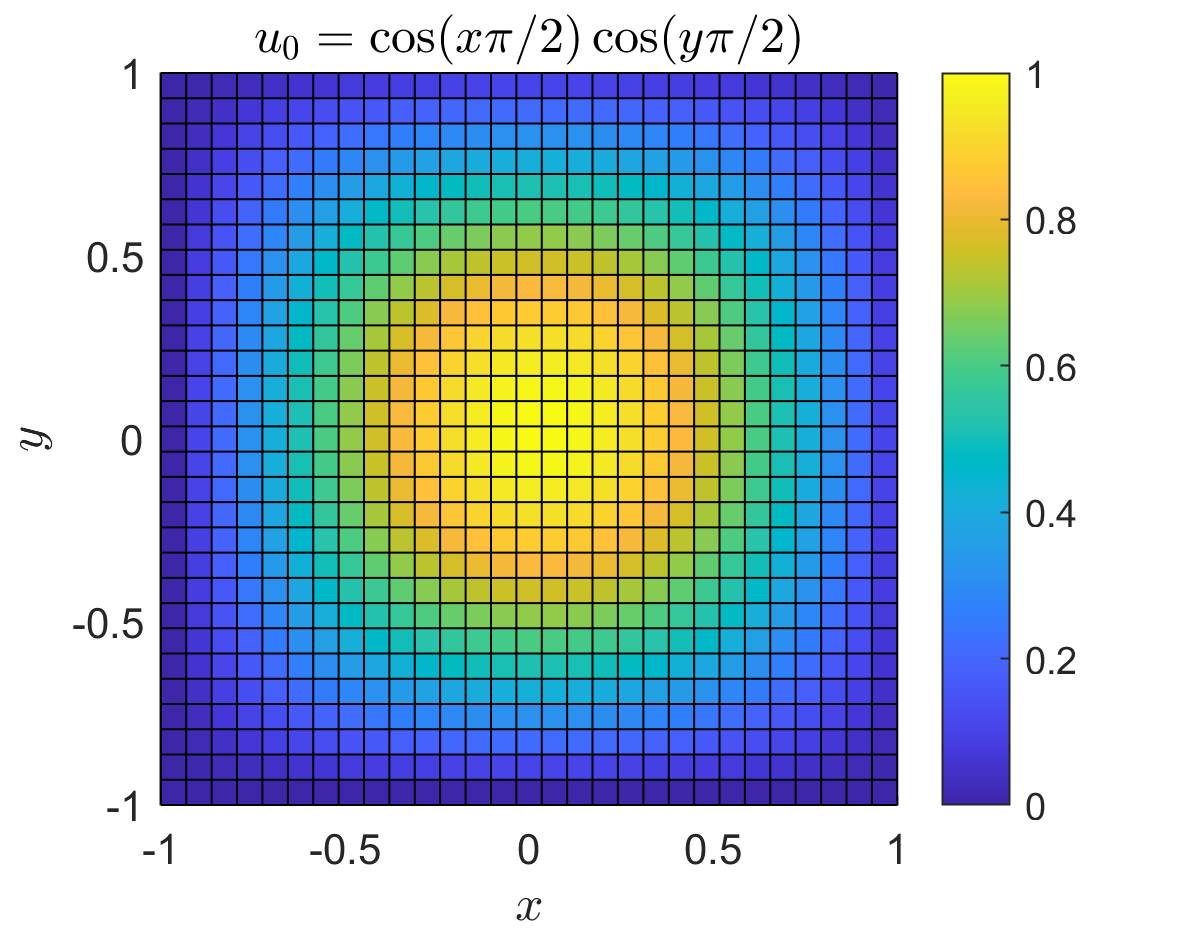}\\
\includegraphics[width=0.45\linewidth]{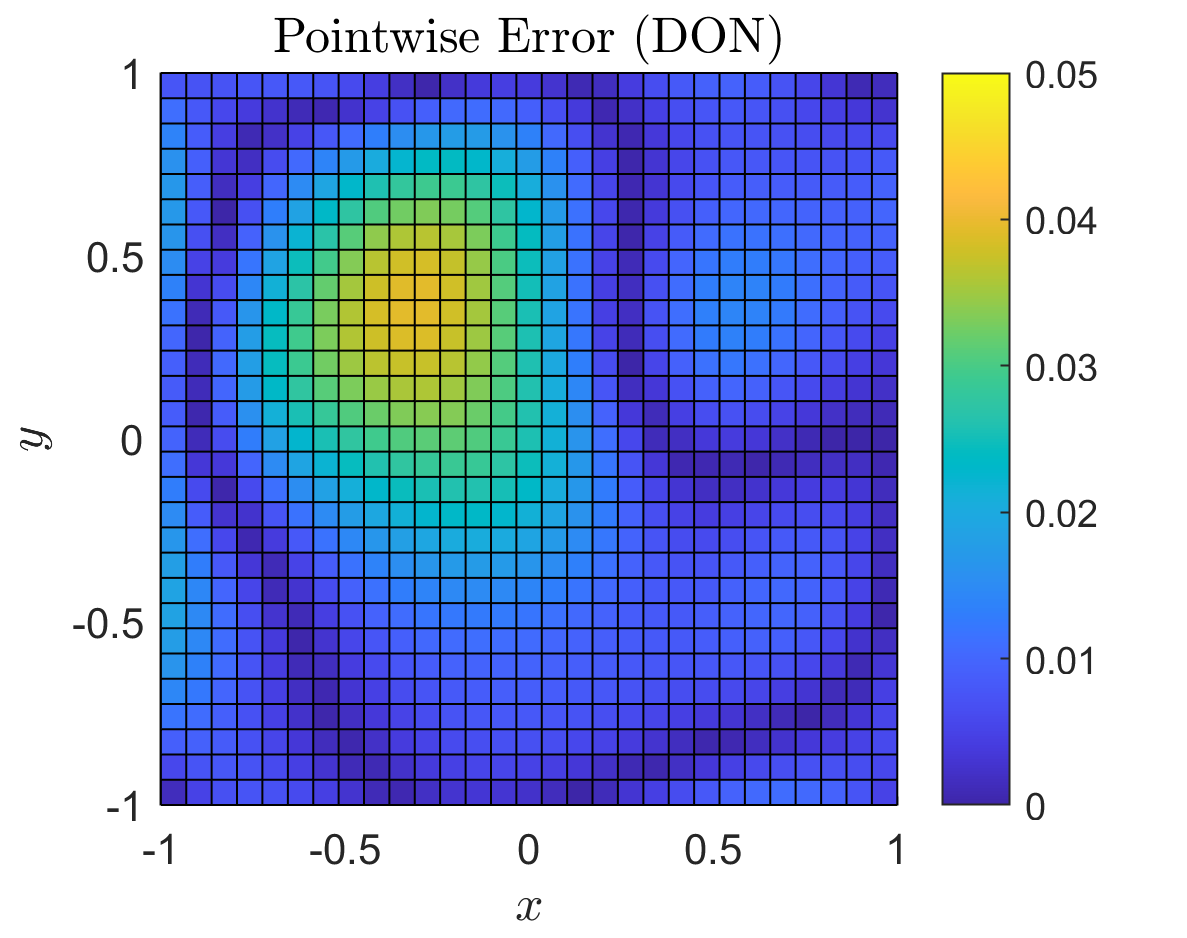}
\includegraphics[width=0.45\linewidth]{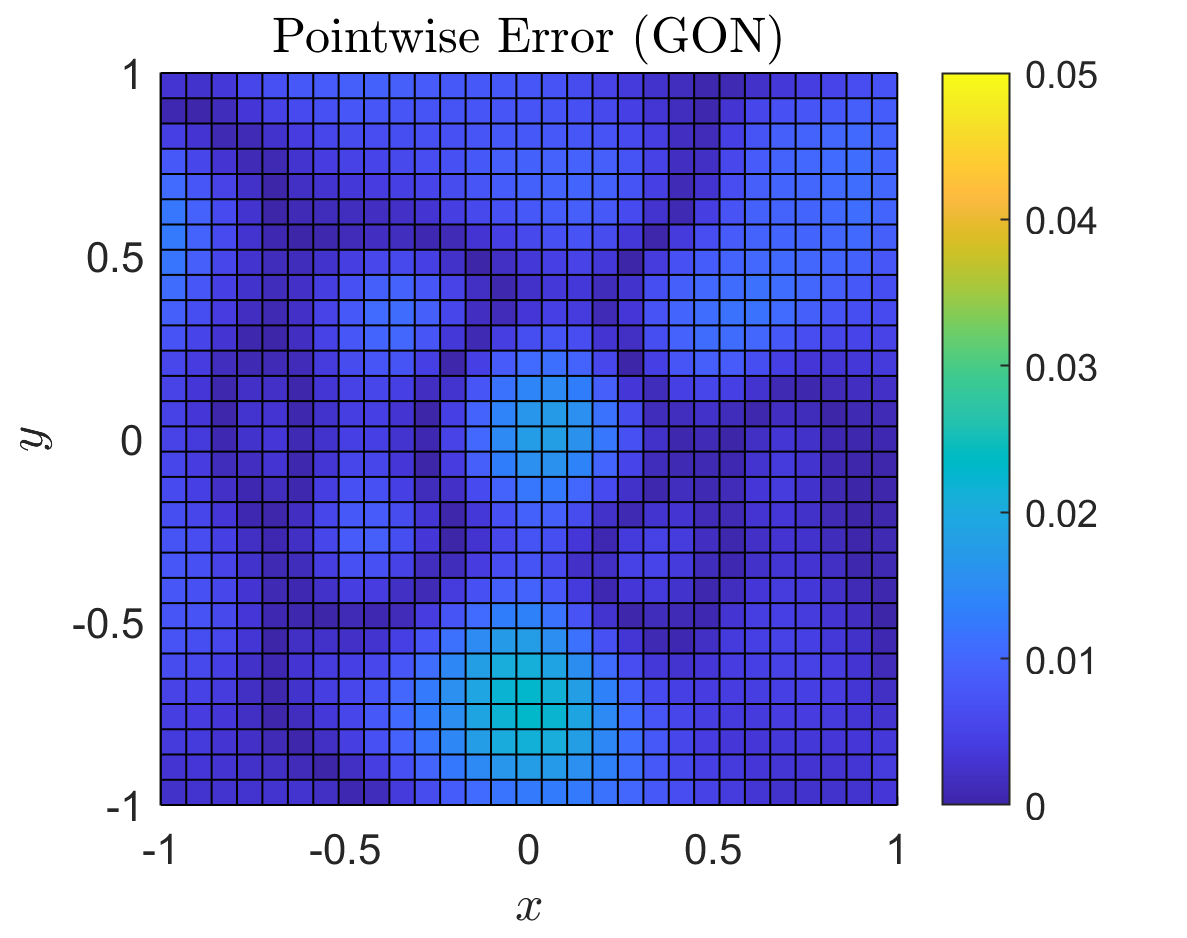}
\caption{Example of Section~\ref{subsec:exp2d}: (left) Initial condition; (middle) Pointwise error at \(t=1.5\) using DeepONets; (right) Pointwise error at \(t=1.5\) using GreenONets.}
\label{fig:func_err2d}
\end{figure}

\begin{figure}[th]
\centering
\includegraphics[width=0.45\linewidth]{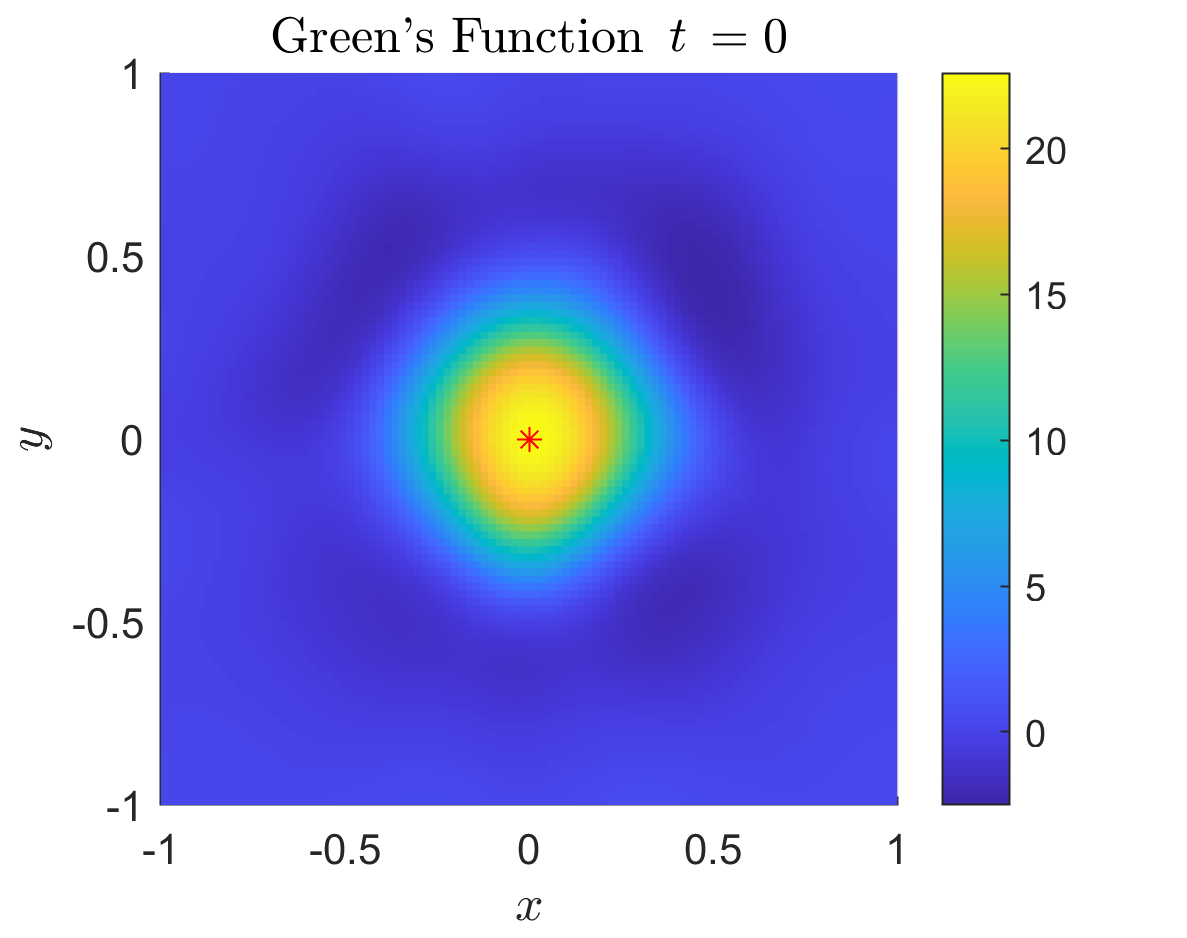}
\includegraphics[width=0.45\linewidth]{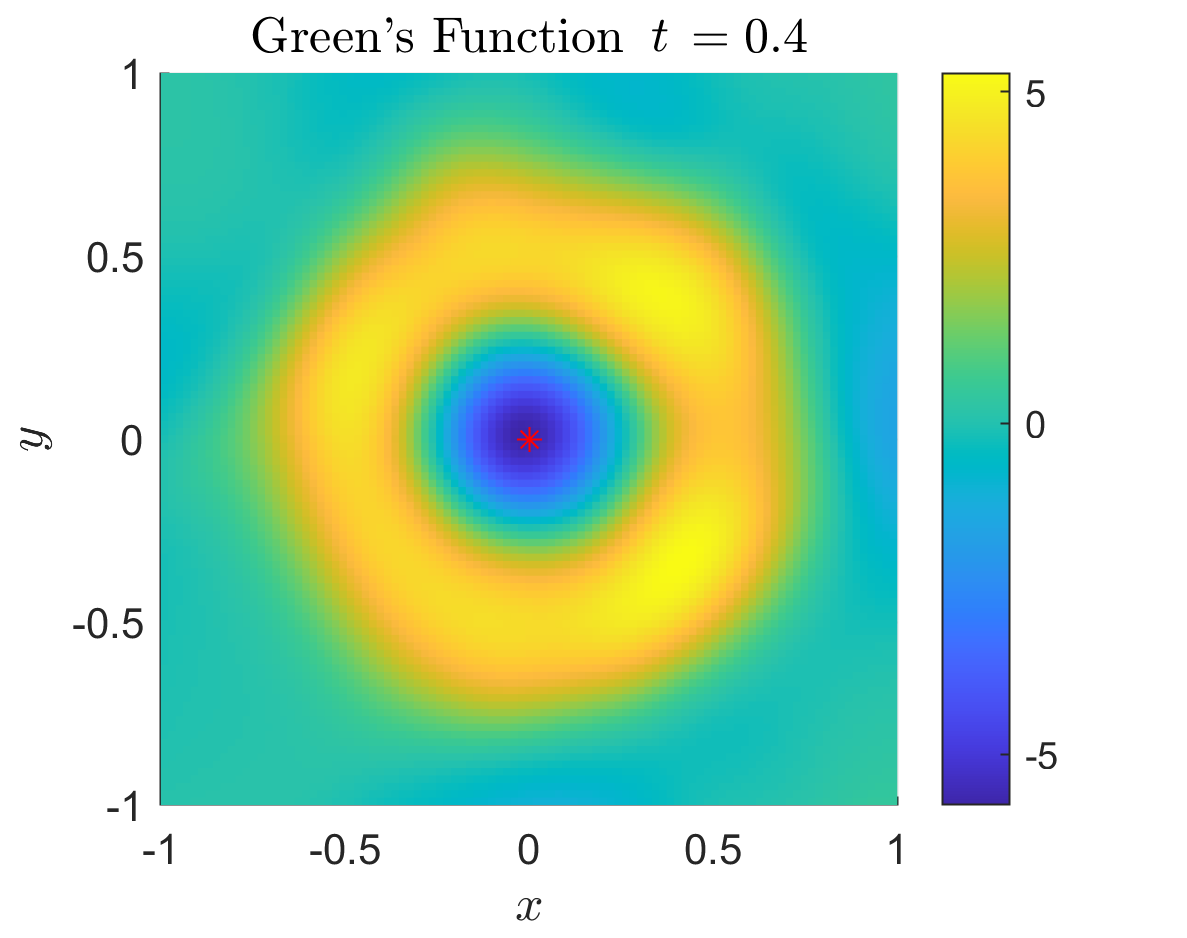}
\caption{Example of Section~\ref{subsec:exp2d}: The approximated Green's function \(G(x,y,t,x_i,y_i) \) for \((x_i,y_i) = (0,0)\) at \(t = 0\) and \(t = 0.4\). The red dot denotes the point \((x_i,y_i) =(0,0)\).}
\label{fig:Greens_2D}
\end{figure}

Similarly to the one-dimensional case, the loss functions have decreased faster for the GreenONets after 400 epochs, see Figure~\ref{fig:loss2d}. We compare in Figure~\ref{fig:func_err2d} the error in the solutions obtained with DeepONets and GreenONets. In this figure, we show the initial condition \(u_0 =\cos(x\pi/2)\cos(y\pi/2)\), with which we test our networks, and the pointwise errors at \(t = 1.5\) using DeepONets and GreenONets. We observe that the maximum pointwise error remains consistently smaller when using GreenONets. Moreover, we plot in Figure~\ref{fig:Greens_2D} the approximated Green's function \(G(x,y,t,x_i,y_i) \) for \((x_i,y_i) = (0,0)\) at \(t = 0\) and \(t = 0.4\). The Green's function initially peaks at the origin and then radially propagates through the domain, as expected. 

\section{Conclusions}
\label{sec:conclusion}

In this work, we have introduced the Green operator networks, that approximate the operator of the wave equation in homogeneous and heterogeneous domains for a family of initial conditions. The GreenONets architecture is inspired by the exact representation of solution of the wave equation in terms of the Green's function in unbounded domains. This architecture yields better results, when approximating the wave operator for homogeneous and heterogeneous domains in one and two dimensions, when compared to DeepONets. The increased performance is attributed to the fact that the GreenONets architecture is better suited to this type of problems, but we also recognize that the DeepONets architecture is more general. We have in particular showed that the loss functions associated with the GreenONets always converge with fewer epochs. The numerical results also highlighted the fact that the pointwise errors are generally smaller with GreenONets and that the solutions generalize better when tested on initial conditions with frequencies higher than that of the training set. Finally, we have observed that the errors in the GreenONets solutions remained localized around the peak amplitudes while the errors with DeepONets had the tendency to spread within the domain. Plans for future works will focus on testing the GreenONets for the two- and three-dimensional wave equation involving materials with various heterogeneous properties and to extend the methodology to other model problems. 
	
\vspace{20pt}
\noindent \textbf{Acknowledgements.} 
SP and ML are grateful for the support from the Natural Sciences and Engineering Research Council of Canada (NSERC) Discovery Grants [grant numbers RGPIN-2019-7154, PGPIN-2018-06592]. This research was also partially supported by an NSERC Collaborative Research and Development Grant [grant number RDCPJ 522310-17] with the Institut de Recherche en \'Electricit\'e du Qu\'ebec and Prompt. ZA and SP are thankful to the Laboratoire de M\'ecanique et d’Acoustique UMR 7031, in Marseille, France, for hosting them. This work received support from the French government under the France 2030 investment plan, as part of the Initiative d'Excellence d'Aix-Marseille Université - A*MIDEX -  AMX-19-IET-010. Finally, the authors also thank J.-P. Ampuero, for having shared SEMLAB~~\cite{Matlabcode} with the community.
	
\bibliographystyle{abbrv}
\bibliography{bibliography}
	
\end{document}